\documentclass{article}
\usepackage{amssymb}
\title{\bf GENERALIZED PASCAL TRIANGLES AND T$\rm \ddot{\bf O}$EPLITZ MATRICES\thanks{This research was in part supported by a
grant from IPM (No. 85200038).}}
\author{{\sc A. R. Moghaddamfar} and {\sc S. M. H. Pooya}\\[0.2cm]
{\em Department of Mathematics, Faculty of Science,}\\
{\em  K. N. Toosi University of Technology},\\
{\em  P. O. Box $16315-1618$, Tehran, Iran}\\ and \\
{\em Institute for Studies in
Theoretical Physics and Mathematics (IPM)}\\
{\em E-mail address}:\\ {\tt moghadam@kntu.ac.ir} \  and \ {\tt
moghadam@mail.ipm.ir}. }

\newenvironment{proof}{\noindent {\em {Proof}}.}{$\square$
\medskip}

\newtheorem{corollary}{Corollary}
\newtheorem{ex}{Example}
\newtheorem{theorem}{Theorem}
\newtheorem{proposition}{Proposition}

\newtheorem{lm}{Lemma}

\begin{document}
\maketitle
\begin{abstract}
\noindent The purpose of this article is to study determinants of
matrices which are known as generalized Pascal triangles (see
\cite{bacher}). We present a factorization by expressing such a
matrix as a product of a unipotent lower triangular matrix, a
T$\rm \ddot{o}$eplitz matrix and a unipotent upper triangular
matrix. The determinant of a generalized Pascal matrix equals
thus the determinant of a T$\rm \ddot{o}$eplitz matrix. This
equality allows us to evaluate a few determinants of generalized
Pascal matrices associated to certain sequences. In particular, we
obtain families of quasi-Pascal matrices whose principal minors
generate any arbitrary linear subsequences $F(nr+s)$ or $L(nr+s)$,
$(n=1, 2, 3, \ldots)$ of Fibonacci or Lucas sequence.
\end{abstract}
\renewcommand{\baselinestretch}{1.1}
\def\thefootnote{ \ }

\section{Introduction}
Let $P(\infty)$ be the infinite symmetric matrix with entries
$P_{i,j}={i+j\choose i}$ for $i,j\geq 0$. The matrix $P(\infty)$
is hence the famous Pascal triangle yielding the binomial
coefficients. The entries of $P(\infty)$ satisfy the recurrence
relation $P_{i,j}=P_{i-1,j}+P_{i,j-1}$. Indeed, this matrix has
the following form:
\begin{equation}\label{13870907} P(\infty)=\left(
\begin{array}{cccccc}
1 & 1& 1& 1& 1 & \ldots\\
1 & 2 & 3& 4& 5 &  \ldots\\
1 & 3 &6 & 10& 15&  \ldots \\
1 & 4 &10 &20 & 35 &  \ldots \\
1 & 5 &15 &35 & 70 &   \ldots \\
\vdots & \vdots & \vdots & \vdots & \vdots & \ddots \\
\end{array} \right)
\end{equation}
One can easily verify that (see \cite{BC, edelman}):
\begin{equation}\label{2008-e1}
P(\infty)=L(\infty)\cdot L(\infty)^t,
\end{equation}
 where $L$ is the infinite unipotent lower
triangular matrix \begin{equation}\label{3000} L(\infty)=\left(
\begin{array}{ccccccc}
1 & & & & &   \\
1 & 1 & & & &  \\
1 & 2 &1 & & &  \\
1 & 3 &3 &1 & &  \\
1 & 4 &6 &4 & 1&  \\
\vdots & \vdots &  \vdots & \vdots & \vdots & \ddots \\
\end{array} \right)
\end{equation}
with entries $L_{i,j}={i \choose j}$.

To introduce the result, we first present some notation and
definitions. We recall that a matrix $T(\infty)=(t_{i,j})_{i,j\geq
0}$ is said to be T$\rm \ddot{o}$eplitz if $t_{i,j}=t_{k,l}$
whenever $i-j=k-l$. Let $\alpha=(\alpha_i)_{i\geq 0}$ and
$\beta=(\beta_i)_{i\geq 0}$ be two sequences with
$\alpha_0=\beta_0$. We shall denote by $T_{\alpha,
\beta}(\infty)=(t_{i,j})_{i,j\geq 0}$ the T$\rm \ddot{o}$eplitz
matrix with $\alpha_i=t_{i,1}$ and $\beta_j=t_{1,j}$. We also
denote by $T_{\alpha, \beta}(n)$ the submatrix of $T_{\alpha,
\beta}(\infty)$ consisting of the elements in its first $n$ rows
and columns.

We come now back to the Eq. (\ref{2008-e1}). In fact, one can
rewrite it as follows:
$$P(\infty)=L(\infty)\cdot I\cdot L(\infty)^t,$$ where matrix $I$ (identity matrix) is a particular case of a
T$\rm \ddot{o}$eplitz matrix.

In \cite{bacher}, Bacher considers determinants of matrices
generalizing the Pascal triangle $P(\infty)$. He introduces
generalized Pascal triangles as follows. Let
$\alpha=(\alpha_i)_{i\geq 0}$ and $\beta=(\beta_i)_{i\geq 0}$ be
two sequences starting with a common first term
$\alpha_0=\beta_0=\gamma$. Then, the generalized Pascal triangle
associated to $\alpha$ and $\beta$, is the infinite matrix
$P_{\alpha,\beta}(\infty)=(P_{i,j})_{i,j\geq 0}$ with entries
$P_{i,0}=\alpha_i$, $P_{0,i}=\beta_i$  and
$$P_{i,j}=P_{i-1,j}+P_{i,j-1}, \ \ {\rm for} \ \  i,j\geq 1.$$
We denote by $P_{\alpha,\beta}(n)$ the finite submatrix of
$P_{\alpha,\beta}(\infty)$ with entries $P_{i,j}$, $0\leq i,j\leq
n-1$.
%Note that $P_{\beta, \alpha}(n)=P_{\alpha,\beta}(n)^t$, and
%so
%$$\det(P_{\alpha,\beta}(n))=\det(P_{\beta,\alpha}(n)).$$
An explicit formula for entry $P_{i,j}$ of $P_{\alpha,\beta}(n)$
is also given by the following formula (see \cite{bacher}):
$$
P_{i,j}=\gamma{i+j\choose
j}+\Big(\sum_{s=1}^{i}(\alpha_{s}-\alpha_{s-1}){i+j-s\choose
j}\Big) +\Big(\sum_{t=1}^{j}(\beta_{t}-\beta_{t-1}){i+j-t\choose
i}\Big).
$$

For arbitrary sequence $\alpha=(\alpha_i)_{i\geq 0}$, we define
the sequences $\hat{\alpha}=(\hat{\alpha}_{i})_{i\geq 0}$ and
$\check{{\alpha}}=(\check{\alpha}_{i})_{i\geq 0}$ as follows:
\begin{equation}\label{2008new}
\hat{\alpha}_{i}=\sum_{k=0}^{i}(-1)^{i+k}{i\choose k}\alpha_k  \
\ \ \mbox{and} \ \ \ \check{\alpha}_{i}=\sum_{k=0}^{i}{i\choose
k}\alpha_k.
\end{equation}

With these definitions we can now state our main result. Indeed,
the purpose of this article is to obtain a {\em factorization} of
the generalized Pascal triangle $P_{\alpha,\beta}(n)$ associated
to the arbitrary sequences $\alpha$ and $\beta$, as a product of
a unipotent lower triangular matrix $L(n)$, a T$\rm
\ddot{o}$eplitz matrix $T_{\hat{\alpha}, \hat{\beta}}(n)$ and a
unipotent upper triangular matrix $U(n)$ (see Theorem \ref{th1}),
that is
$$P_{\alpha,\beta}(n)=L(n)\cdot T_{\hat{\alpha},
\hat{\beta}}(n)\cdot
 U(n).$$
Similarly, we show that
$$T_{\alpha,\beta}(n)=L(n)^{-1}\cdot P_{\check{\alpha},
\check{\beta}}(n)\cdot
 U(n)^{-1}.$$
In fact, we obtain a {\em connection} between generalized Pascal
triangles and T$\rm \ddot{o}$eplitz matrices. In view of these
factorizations, we can easily see that
$$\det(P_{\alpha,\beta}(n))=\det(T_{\hat{\alpha}, \hat{\beta}}(n)).$$ Finally, we present {\em
several applications} of Theorem \ref{th1} to some other
determinant evaluations.

We conclude the introduction with notation and terminology to be
used throughout the article. By $\lfloor x\rfloor$ we denote the
integer part of $x$, i. e., the greatest integer that is less than
or equal to $x$. We also denote by $\lceil x\rceil$ the smallest
integer greater than or equal to $x$. Given a matrix $A$, we
denote by ${\rm R}_i(A)$ and ${\rm C}_j(A)$ the row $i$ and the
column $j$ of $A$, respectively. We use the notation $A^t$ for the
transpose of $A$. Also, we denote by
$$A^{j_1,j_2,\ldots,j_k}_{i_1,i_2,\ldots,i_l}$$ the submatrix of
$A$ obtained by erasing rows $i_1,i_2,\ldots,i_l$ and columns
$j_1,j_2,\ldots,j_k$. In general, an $n\times n$  matrix of the
following form:
$$\left [\begin{array}{l|c} A & B \\
\hline C & P_{\alpha, \beta}(n-k) \end{array} \right] \Big({\rm
resp.} \left [\begin{array}{l|c} A & B \\
\hline C & T_{\alpha, \beta}(n-k) \end{array} \right]\Big)$$ where
$A$, $B$ and $C$ are arbitrary matrices of order ${k\times k}$,
${k\times (n-k)}$ and ${(n-k)\times k}$, respectively, is called
a {\em quasi-Pascal} (resp. {\em quasi-T$\rm \ddot{o}$eplitz})
matrix.

Throughout this article we assume that:
\begin{center}
\begin{tabular}{ll}
$\mathcal{F}=(\mathcal{F}_i)_{i\geq
0}=(0, 1, 1, 2, 3, 5, 8, \ldots, \mathcal{F}_i=F(i), \ldots)$ & (Fibonacci  numbers),\\[0.15cm]
$\mathcal{F}^\ast=(\mathcal{F}_i)_{i\geq
1}=(1, 1, 2, 3, 5, 8, \ldots, \mathcal{F}_i=F(i), \ldots)$ & (Fibonacci  numbers $\neq 0$),\\[0.15cm]

$\mathcal{L}=(\mathcal{L}_i)_{i\geq
0}=(2, 1, 3, 4, 7, 11, 18, \ldots, \mathcal{L}_i=L(i), \ldots)$ & (Lucas  numbers),\\[0.15cm]
$\mathcal{C}=(\mathcal{C}_i)_{i\geq
0}=(1, 1, 2, 5, 14, 42, 132, \ldots, \mathcal{C}_i=C(i), \ldots)$ & (Catalan  numbers),\\[0.15cm]
$\mathcal{I}=(\mathcal{I}_i)_{i\geq
0}=(0!, 1!, 2!, 3!, 4!, 5!, 6!, \ldots, \mathcal{I}_i=i!, \ldots)$ & \\[0.15cm]
$\mathcal{I}^\ast=(\mathcal{I}_i)_{i\geq
1}=(1!, 2!, 3!, 4!, 5!, 6!, \ldots, \mathcal{I}_i=i!, \ldots)$ & \\[0.15cm]
\end{tabular}
\end{center}
The paper is organized as follows: In Section 2, we derive some
preparatory results. In Section 3, we prove the main result
(Theorem 1) and in Section 4, we present some applications of
main theorem.

\section{Preliminary Results}
As we mentioned in Introduction, if $\alpha=(\alpha_i)_{i\geq 0}$
is an arbitrary sequence, then we define the sequences
$\hat{\alpha}=(\hat{\alpha}_{i})_{i\geq 0}$ and
$\check{{\alpha}}=(\check{\alpha}_{i})_{i\geq 0}$ as in Eq.
(\ref{2008new}). For some certain sequences $\alpha$ the
associated sequences $\hat{\alpha}$ and $\check{{\alpha}}$ seem
also to be of interest since they have appeared elsewhere. In
Tables 1 and 2, we have presented some sequences $\alpha$ and the
associated sequences $\hat{\alpha}$ and $\check{{\alpha}}$.
\begin{center}
{\bf Table 1.} Some well-known sequences $\alpha$ and associated
sequences $\check{{\alpha}}$.
\end{center}
$$\begin{tabular}{|lc|l|}
\hline $\alpha$ &  Reference  & $\check{{\alpha}}$  \\
\hline  $(0, 1, -1, 2, -3, 5, -8, \ldots )$&  ${\rm \bf A}039834$
in
\cite{IS} & $\mathcal{F}$   \\
\hline  $(1, 0, 1, -1, 2, -3, 5, -8, \ldots )$&  ${\rm \bf
A}039834$ in
\cite{IS} & $\mathcal{F}^\ast$  \\
\hline  $(2, -1, 3, -4, 7, -11, 18, \ldots )$&  ${\rm \bf
A}061084$ in
\cite{IS} & $\mathcal{L}$  \\
\hline  $(1, 0, 1, 1, 3, 6, 15, \ldots )$&   ${\rm \bf A}005043$
in
\cite{IS} & $\mathcal{C}$   \\
\hline  $(1, 0, 1, 2, 9, 44, 265, \ldots )$&  ${\rm \bf A}000166$
in
\cite{IS}  & $\mathcal{I}$  \\
\hline  $(1, 1, 3, 11, 53, 309, 2119, \ldots )$&  ${\rm \bf
A}000255$ in \cite{IS} & $\mathcal{I}^\ast$   \\ \hline
\end{tabular}
$$
\begin{center}
{\bf Table 2.} Some well-known sequences $\alpha$ and associated
sequences $\hat{\alpha}$.
\end{center}
$$
\begin{tabular}{|lc|l|}
\hline  $\alpha$ &  Reference  & $\hat{\alpha}$   \\
\hline
$(0, 1, 3, 8, 21, 55, 144, \ldots)$ &  ${\rm \bf A}001906$ in \cite{IS} & $\mathcal{F}$ \\
$(1, 2, 5, 13, 34, 89, 233, \ldots)$ &  ${\rm \bf A}001519$ in \cite{IS} & $\mathcal{F}^\ast$ \\
$(2, 3, 7, 18, 47, 123, 322, \ldots)$ &  ${\rm \bf A}005248$ in \cite{IS} & $\mathcal{L}$ \\
$(1, 2, 5, 15, 51, 188, 731, \ldots)$ &  ${\rm \bf A}007317$ in \cite{IS} & $\mathcal{C}$ \\
$(1, 2, 5, 16, 65, 326, 1957, \ldots)$ &  ${\rm \bf A}000522$ in \cite{IS} &  $\mathcal{I}$ \\
$(1, 3, 11, 49, 261, 1631, \ldots)$ &  ${\rm \bf A}001339$ in \cite{IS} & $\mathcal{I}^\ast$  \\
\hline
\end{tabular}
$$

As an another nice example, consider $\alpha={\rm
R}_i(L(\infty))$ where $L$ is introduced as in Eq. (\ref{3000}).
Then, we have
$$\check{\alpha}={\rm R}_i(P(\infty)).$$
\begin{lm}\label{2008-lm1}
Let $\alpha$ be a sequence. Then we have
$\hat{\check{\alpha}}={\check{\hat \alpha}}=\alpha.$
\end{lm}
\begin{proof} Suppose $\alpha=(\alpha_i)_{i\geq 0}$ and $\hat{\check{\alpha}}=(\hat{\check{\alpha}}_i)_{i\geq
0}$. Then, we have
$$
\begin{tabular}{lll}
$\hat{\check{\alpha}}_i$ &= & $\sum_{k=0}^{i}(-1)^{i+k}{i\choose k}\check{\alpha}_k$\\[0.2cm]
&= &$\sum_{k=0}^{i}(-1)^{i+k}{i\choose k}\sum_{s=0}^{k}{k\choose s}\alpha_s $\\[0.2cm]
&= &$\sum_{k=0}^{i}\sum_{s=0}^{k}(-1)^{i+k}{i\choose k}{k\choose s}\alpha_s $\\[0.2cm]
&= &$(-1)^{i}\sum_{l=0}^{i}\alpha_{l}\sum_{h=l}^{i}(-1)^{h}{i\choose h}{h\choose l}$\\[0.2cm]
&= &$(-1)^{i}\sum_{l=0}^{i}\alpha_{l}\sum_{h=l}^{i}(-1)^{h}{i\choose l}{i-l \choose h-l}$\\[0.2cm]
&= &$(-1)^{i}\sum_{l=0}^{i}\alpha_{l}{i\choose l}\sum_{h=l}^{i}(-1)^{h}{i-l \choose h-l}$.\\[0.2cm]
\end{tabular}
$$
But, if $l<i$, then we have $\sum_{h=l}^{i}(-1)^{h}{i-l \choose
h-l}=0$. Therefore, we obtain
$$
\hat{\check{\alpha}}_i= (-1)^{i}\sum_{l=i}^{i}\alpha_{l}{i\choose
l}\sum_{h=l}^{i}(-1)^{h}{i-l \choose h-l}=\alpha_{i},$$ and hence
$\hat{\check{\alpha}}=\alpha.$

The proof of second part is similar to the previous case.
\end{proof}

\begin{lm}\label{ele-1}
Let $i, j$ be positive integers. Then we have
\[ \sum_{k=0}^{i-j}(-1)^{k}{i\choose k+j}{k+j \choose j}
=\left \{ \begin{array}{ll} 0 &  \ \mbox{if} \ \ \  i\neq j,\\[0.2cm] 1 &
\ \mbox{if} \ \ \  i=j.
\end{array} \right. \]
\end{lm}
The {\em proof} follows from the easy identity $${i\choose
k+j}{k+j\choose j}={i\choose j}{i-j\choose k}.$$

The following Lemma is a special case of a general result due to
Krattenthaler (see Theorem 1 in \cite{kkk}).
\begin{lm}\label{ele1387-18}
Let $\alpha=(\alpha_i)_{i\geq 0}$ and $\beta=(\beta_j)_{j\geq 0}$
be two geometric sequences with $\alpha_i=\rho^{i}$ and
$\beta_j=\sigma^j$. Then, we have
$$\det(P_{\alpha, \beta}(n))=(\rho+\sigma-\rho\sigma)^{n-1}.$$
\end{lm}

\section{Main Result}
Now, we are in the position to state and prove the main result of
this article.
\begin{theorem}\label{th1}
Let $\alpha=(\alpha_i)_{i\geq 0}$ and $\beta=(\beta_i)_{i\geq 0}$
be two sequences starting with a common first term
$\alpha_0=\beta_0=\gamma$. Then we have
\begin{equation}\label{e}
P_{\alpha,\beta}(n)=L(n)\cdot T_{\hat{\alpha},
\hat{\beta}}(n)\cdot
 U(n),
\end{equation}
and
\begin{equation}\label{2008-e}
T_{\alpha,\beta}(n)=L(n)^{-1}\cdot P_{\check{\alpha},
\check{\beta}}(n)\cdot
 U(n)^{-1},
\end{equation}
where $L(n)=(L_{i,j})_{0\leq i,j<n}$ is a lower triangular matrix
with
$$L_{i,j}=\left \{ \begin{array}{ccc}
0 & \mbox{if} & i<j\\
{i\choose j} & \mbox{if} & i\geq j,\\
\end{array} \right.
$$
$U(n)=L(n)^t$. In particular, we have
$\det(P_{\alpha,\beta}(n))=\det(T_{\hat{\alpha},
\hat{\beta}}(n))$.
\end{theorem}

\begin{proof} First, we claim
that $$P_{\alpha,\beta}(n)=L(n)\cdot Q(n),$$ where
$L(n)=(L_{i,j})_{0\leq i,j<n}$ is a lower triangular matrix with
$$L_{i,j}=\left \{ \begin{array}{ccc}
0 & \mbox{if} & i<j\\
{i\choose j} & \mbox{if} & i\geq j,\\
\end{array} \right.
$$
and $Q(n)=(Q_{i,j})_{0\leq i,j<n}$ with $Q_{i,0}=\hat{\alpha}_i$,
$Q_{0,i}=\beta_i$ and
\begin{equation}\label{ee}
Q_{i,j}=Q_{i-1,j-1}+Q_{i,j-1}, \ \ \ \ 1\leq i,j<n.
\end{equation}
For instance, when $n=4$ the matrices $L(4)$ and $Q(4)$ are given
by:
$$L(4)=\left[ \begin{tabular}{llll}
$1$ & $0$ & $0$ & $0$ \\
$1$ & $1$ & $0$ & $0$ \\
1 & $2$ & $1$ & $0$ \\
1 & 3 & $3$ & $1$ \\
\end{tabular} \right ], $$
and
$$
Q(4)=\left[ \begin{array}{cccc}
\gamma & \beta_1 & \beta_2 & \beta_3 \\
-\gamma+\alpha_1 & \alpha_1 & \beta_1+\alpha_1 & \beta_1+\beta_2+\alpha_2 \\
\gamma-2\alpha_1+\alpha_2 & -\alpha_1+\alpha_2 & \alpha_2 & \beta_1+\alpha_1+\alpha_2 \\
-\gamma+3\alpha_1-3\alpha_2+\alpha_3 & \alpha_1-2\alpha_2+\alpha_3 & -\alpha_2+\alpha_3 & \alpha_3 \\
\end{array} \right ].$$
Note that the entries of $L(n)$ satisfying in the following
recurrence
\begin{equation}\label{e2}
L_{i,j}=L_{i-1,j-1}+L_{i-1,j}, \ \ \ 1\leq i,j<n.
\end{equation}
For the proof of the claimed factorization we compute the
$(i,j)$-th entry of $L(n)\cdot Q(n)$, that is
$$(L(n)\cdot Q(n))_{i,j}=\sum_{k=1}^nL_{i,k}Q_{k,j}.$$ In fact, it suffices to show that
\begin{center}
\begin{tabular}{l}
${\rm R}_0(L(n)\cdot Q(n))={\rm
R}_0(P_{\alpha,\beta}(n)),$\\[0.2cm]
${\rm C}_0(L(n)\cdot Q(n))={\rm C}_0(P_{\alpha,\beta}(n))$
\end{tabular}
\end{center}
 and
\begin{equation}\label{e3}
(L(n)\cdot Q(n))_{i,j}=(L(n)\cdot Q(n))_{i,j-1}+(L(n)\cdot
Q(n))_{i-1,j},
\end{equation}
for $1\leq i, j<n$.

First, suppose that $i=0$. Then
$$(L(n)\cdot Q(n))_{0,j}=\sum_{k=0}^{n-1}L_{0,k}Q_{k,j}=L_{0,0}Q_{0,j}=\beta_j,$$
and so ${\rm R}_0(L(n)\cdot Q(n))={\rm
R}_0(P_{\alpha,\beta}(n))=({\beta}_0, \ {\beta}_1, \ \ldots, \
{\beta}_{n-1})$.

Next, suppose that $i\geq 1$ and $j=0$. In this case, we have
$$
\begin{tabular}{lll}
$(L(n)\cdot Q(n))_{i,0}$ & $=$ & $\sum_{k=0}^{n-1}L_{i,k}Q_{k,0}$
\\[0.2cm]
& $=$ & $\sum_{k=0}^{i}\Big\{{i\choose k}
\sum_{l=0}^{k}(-1)^{l+k}{k\choose l}\beta_{l}\Big \}$\\[0.2cm]
& $=$ &
$\sum_{j=0}^{i}\beta_{j}\Big\{\sum_{t=0}^{i-j}(-1)^t{i \choose t+j}{t+j \choose j}\Big \}$\\[0.2cm]
& $=$ & $\beta_i$,  \ \ \ (by Lemma \ref{ele-1})
\end{tabular}
$$ and so ${\rm C}_0(L(n)\cdot Q(n))={\rm
C}_0(P_{\alpha, \beta}(n))=(\alpha_0, \ \alpha_1, \ \ldots, \
\alpha_{n-1})$.

Finally, we must establish Eq. (\ref{e3}). Therefore, we assume
that $1\leq i, j<n$. In this case, we have
$$
\begin{tabular}{lll} $(L(n)\cdot Q(n))_{i,j}$ & =&
$\sum_{k=0}^{n-1}L_{i,k}Q_{k,j}$\\[0.2cm]
& = & $L_{i,0}Q_{0,j}+\sum_{k=1}^{n-1}L_{i,k}Q_{k,j}$\\[0.2cm]
& = &
$L_{i,0}Q_{0,j}+\sum_{k=1}^{n-1}L_{i,k}(Q_{k-1,j-1}+Q_{k,j-1})$ \ \ (by Eq. (\ref{ee}) )\\[0.2cm]
& =& $L_{i,0}Q_{0,j}+\sum_{k=1}^{n-1}L_{i,k}Q_{k-1,j-1}+\sum_{k=1}^{n-1}L_{i,k}Q_{k,j-1}$\\[0.2cm]
& = &
$L_{i,0}Q_{0,j}+\sum_{k=1}^{n-1}(L_{i-1,k-1}+L_{i-1,k})Q_{k-1,j-1}$ \\[0.2cm]
 & & $+\sum_{k=0}^{n-1}L_{i,k}Q_{k,j-1}-L_{i,0}Q_{0,j-1}$ (by Eq. (\ref{e2}))\\[0.2cm] & = &
$L_{i,0}Q_{0,j}+\sum_{k=1}^{n-1}L_{i-1,k-1}Q_{k-1,j-1}+\sum_{k=1}^{n-1}L_{i-1,k}(Q_{k,j}-Q_{k,j-1})$\\[0.2cm]
& & $+(L(n)\cdot Q(n))_{i,j-1}-L_{i,0}Q_{0,j-1}$ \ \ \ (by Eq.
(\ref{ee}))
\\[0.2cm]
& = &
$L_{i,0}Q_{0,j}+\sum_{k=0}^{n-2}L_{i-1,k}Q_{k,j-1}+\sum_{k=1}^{n-1}L_{i-1,k}Q_{k,j}$ \\[0.2cm]
&  &
$-\sum_{k=0}^{n-1}L_{i-1,k}Q_{k,j-1}+L_{i-1,0}Q_{0,j-1}+(L(n)\cdot
Q(n))_{i,j-1}-L_{i,0}Q_{0,j-1}$
\\[0.2cm]
& = &
$L_{i,0}Q_{0,j}+\sum_{k=0}^{n-1}L_{i-1,k}Q_{k,j-1}+\sum_{k=0}^{n-1}L_{i-1,k}Q_{k,j}-L_{i-1,0}Q_{0,j}$ \\[0.2cm]
&  &
$-\sum_{k=0}^{n-1}L_{i-1,k}Q_{k,j-1}+L_{i-1,0}Q_{0,j-1}+(L(n)\cdot
Q(n))_{i,j-1}-L_{i,0}Q_{0,j-1}$
\\[0.2cm]
& & \ \ \ (note that $L_{i-1,n-1}=0$)\\[0.2cm]
& = & $(L_{i,0}-L_{i-1,0})Q_{0,j}+(L(n)\cdot
Q(n))_{i-1,j}+(L_{i-1,0}-L_{i,0})Q_{0,j-1}$\\[0.2cm]
& & $+(L(n)\cdot
Q(n))_{i,j-1}$ \\[0.2cm]

& = &
$(L(n)\cdot Q(n))_{i-1,j}+(L(n)\cdot Q(n))_{i,j-1}$,  \ \ (note that $L_{i,0}=L_{i-1,0}=1$)\\[0.2cm]
\end{tabular}
$$
which is Eq. (\ref{e3}).

Next, we claim that $$Q(n)=T_{\hat{\alpha}, \hat{\beta}}(n)\cdot
U(n),$$ where $U(n)=L(n)^t$ and $T_{\alpha,
\beta}(n)=(T_{i,j})_{0\leq i,j<n}$ with
$T_{i,0}=\hat{\alpha}_{i}$, $T_{0,j}=\hat{\beta}_{j}$, and
$$T_{i,j}=T_{i-1,j-1}, \ \ \ 1\leq i,j<n.$$
Note that, we have
\begin{equation}\label{e4}
U_{i,j}=U_{i-1,j-1}+U_{i,j-1}, \ \ \ 1\leq i,j<n.
\end{equation}
For instance, when $n=4$ the matrices $T_{\hat{\alpha},
\hat{\beta}}(4)$ and $U(4)$ are given by:
$$
T_{\hat{\alpha}, \hat{\beta}}(4)=\left[ \begin{array}{cccc}
\gamma & -\gamma+\beta_1 & \gamma-2\beta_1+\beta_2 & -\gamma+3\beta_1-3\beta_2+\beta_3\\
-\gamma+\alpha_1 & \gamma& -\gamma+\beta_1  & \gamma-2\beta_1+\beta_2 \\
\gamma-2\alpha_1+\alpha_2 & -\gamma+\alpha_1& \gamma & -\gamma+\beta_1  \\
-\gamma+3\alpha_1-3\alpha_2+\alpha_3 & \gamma-2\alpha_1+\alpha_2 & -\gamma+\alpha_1& \gamma \\
\end{array} \right ],$$
and
$$U(4)=\left[ \begin{tabular}{llll}
$1$ & $1$ & $1$ & $1$ \\
$0$ & $1$ & $2$ & $3$ \\
$0$ & $0$ & $1$ & $3$ \\
$0$ & $0$ & $0$ & $1$ \\
\end{tabular} \right ]. $$

As before, the proof of the claim requires some calculations. If
we have $i=0$, then
$$
\begin{tabular}{lll}
$(T_{\hat{\alpha}, \hat{\beta}}(n)\cdot U(n))_{0,j}$ & $=$ &
$\sum_{k=0}^{n-1}T_{0,k}U_{k,j}$
\\[0.2cm]
& $=$ & $\sum_{k=0}^{j}\Big\{\Big(
\sum_{l=0}^{k}(-1)^{l+k}{k\choose l}\beta_{l}\Big){j\choose k}\Big \}$\\[0.2cm]
& $=$ &
$\sum_{i=0}^{j}\beta_{i}\Big\{\sum_{t=0}^{j-i}(-1)^t{j \choose t+i}{t+i \choose i}\Big \}$\\[0.2cm]
& $=$ & $\beta_j$,  \ \ \ (by Lemma \ref{ele-1})
\end{tabular}
$$
which implies that ${\rm R}_0(T_{\hat{\alpha},
\hat{\beta}}(n)\cdot U(n))={\rm R}_0(Q(n))=({\beta}_0, \
{\beta}_1, \ \ldots, \ {\beta}_{n-1})$. If $j=0$, then we obtain
$$(T_{\hat{\alpha}, \hat{\beta}}(n)\cdot
U(n))_{i,0}=\sum_{k=0}^{n-1}T_{i,k}U_{k,1}=T_{i,0}U_{0,0}=\hat{\alpha_i},
$$
and so ${\rm C}_0(T_{\hat{\alpha}, \hat{\beta}}(n)\cdot U(n))={\rm
C}_0(Q(n))=(\hat{\alpha}_0, \ \hat{\alpha}_1, \ \ldots, \
\hat{\alpha}_{n-1})$. Finally, we assume that $1\leq i,j<n-1$ and
establish Eq. (\ref{ee}). Indeed, by calculations we observe that
$$
\begin{tabular}{lll}
$(T_{\hat{\alpha}, \hat{\beta}}(n)\cdot U(n))_{i,j}$ & = &
$\sum_{k=0}^{n-1}T_{i,k}U_{k,j}
$\\[0.2cm]
& = &
$T_{i,0}U_{0,j}+\sum_{k=1}^{n-1}T_{i,k}U_{k,j}$\\[0.2cm]
& = &  $T_{i,0}U_{0,j}+\sum_{k=1}^{n-1}T_{i,k}(U_{k-1,j-1}+U_{k,j-1})$ \ \ \ (by Eq. (\ref{e4}))\\[0.2cm]
& = &  $T_{i,0}U_{0,j}+\sum_{k=1}^{n-1}T_{i,k}U_{k-1,j-1}+\sum_{k=1}^{n-1}T_{i,k}U_{k,j-1}$\\[0.2cm]
& = &  $T_{i,0}U_{0,j}+\sum_{k=1}^{n-1}T_{i-1,k-1}U_{k-1,j-1}+
\sum_{k=0}^{n-1}T_{i,k}U_{k,j-1}-T_{i,0}U_{0,j-1}$  \\[0.2cm]
& & (note that, $T_{i,k}=T_{i-1,k-1}$)  \\[0.2cm] & = &
$T_{i,0}(U_{0,j}-U_{0,j-1})+\sum_{k=0}^{n-1}T_{i-1,k}U_{k,j-1}+
(T_{\hat{\alpha}, \hat{\beta}}(n)\cdot U(n))_{i,j-1}$\\[0.2cm]
&& (note that $U_{n-1,j-1}=0$)\\[0.2cm]
& = &  $(T_{\hat{\alpha}, \hat{\beta}}(n)\cdot U(n))_{i-1,j-1}+(T_{\hat{\alpha}, \hat{\beta}}(n)\cdot U(n))_{i,j-1}$,\\[0.2cm]
\end{tabular}
$$
which is Eq. (\ref{ee}). The proof of Eq. (\ref{e}) is now
complete.

To prove of Eq. (\ref{2008-e}), we observe that
$$
\begin{tabular}{lll}
$L(n)^{-1}\cdot P_{\check{\alpha}, \check{\beta}}(n)\cdot
U(n)^{-1}$ &= & $L(n)^{-1}\cdot L(n)\cdot T_{\hat{\check{\alpha}},
\hat{\check{\beta}}}(n)\cdot U(n) \cdot U(n)^{-1}$ \ \ \ \ \ (by Eq. (\ref{e}))\\[0.2cm]
&= &$ T_{\hat{\check{\alpha}},
\hat{\check{\beta}}}(n)$\\[0.2cm]
&= &$ T_{\alpha, \beta}(n)$.  \ \ \ \ \ (by Lemma \ref{2008-lm1}) \\[0.2cm]
\end{tabular}
$$
The proof is complete.
\end{proof}

\section{Some Applications}

\subsection{Generalized Pascal Triangle Associated to an Arithmetical or Geometric Sequence }

\begin{corollary}\label{th2}
Let $\alpha=(\alpha_i)_{i\geq 0}$ be an arithmetical sequence with
$\alpha_i=a+id$, and let $\beta=(\beta_i)_{i\geq 0}$ be an
arbitrary sequence with $\beta_0=a$. We set $D(n)=\det(P_{\alpha,
\beta}(n))$. Then we have
$$D(n)=\sum_{k=0}^{n-1}(-d)^{k}\hat{\beta}_kD(n-k-1),$$
with $D(0)=1$.
\end{corollary}
\begin{proof}
By Theorem \ref{th1}, we deduce that $D(n)=\det(T_{\hat{\alpha},
\hat{\beta}}(n))$, where $\hat{\alpha}=(\hat{\alpha}_i)_{i\geq
0}$ with
$$
\begin{tabular}{lll}
$\hat{\alpha}_i$ & $=$ &
$\sum_{k=0}^{i}(-1)^{i+k}{i \choose k}\alpha_k$ \\[0.2cm]
& $=$ &
$\sum_{k=0}^{i}(-1)^{i+k}{i \choose k}\big(\beta_0+kd\big)$\\[0.2cm]
& $=$ &
$\beta_0\sum_{k=0}^{i}(-1)^{i+k}{i \choose k}+d\sum_{k=0}^{i}(-1)^{i+k}{i \choose k}k$\\[0.2cm]
& $=$ & $\left \{ \begin{array}{ccc} \beta_0 & \mbox{if}
& i=0,\\[0.1cm]
d & \mbox{if} & i=1,\\[0.1cm]
0 & \mbox{if} & i>1.\\ \end{array} \right.$\\
\end{tabular}
$$
Now, expanding through the first row of $T_{\hat{\alpha},
\hat{\beta}}(n)$, we obtain the result. \end{proof}

\begin{corollary}\label{cor2}
Let $\alpha=(\alpha_i)_{i\geq 0}$ be an arithmetical sequence with
$\alpha_i=a+id$, and let $\beta=(\beta_{i})_{i\geq 0}$ be an
alternating sequence with $\beta_i=(-1)^{i}a$. Then we have
$$\det(P_{\alpha, \beta}(n))=a(2d+a)^{n-1}.$$
\end{corollary}
\begin{proof}
Let $D(n)=\det(P_{\alpha, \beta}(n))$. By Corollary \ref{th2}, we
have
$$D(n)=\sum_{k=0}^{n-1}(-d)^{k}\hat{\beta}_kD(n-k-1) \ \ \ \ \  (n\geq 1),$$
with $D(0)=1$. An easy calculation shows that
$$
\begin{tabular}{lll}
$\hat{\beta}_k$ & $=$ &
$\sum_{l=0}^{k}(-1)^{k+l}{k \choose l}\beta_l$ \\[0.2cm]
& $=$ &
$\sum_{l=0}^{k}(-1)^{k+l}{k \choose l}(-1)^{l}a$\\[0.2cm]
& $=$ &
$a(-1)^{k}\sum_{l=0}^{k}{k \choose l}$\\[0.2cm]
& $=$ & $a(-2)^{k}.$\\
\end{tabular}
$$
Hence, we have
$$ D(n)=a\sum_{k=0}^{n-1}(2d)^{k}D(n-k-1).$$
Now, replacing $n$ by $n-1$, we obtain
$$ D(n-1)=a\sum_{k=0}^{n-2}(2d)^{k-1}D(n-k-2),$$
and by calculation it follows that
$$
\begin{tabular}{lll}
$D(n)-2dD(n-1)$ & $=$ & $a\Big(\sum_{k=0}^{n-1}(2d)^{k}D(n-k-1)-\sum_{k=0}^{n-2}(2d)^{k-1}D(n-k-2)\Big) $\\[0.2cm]
& $=$ &
$a\Big(\sum_{k=0}^{n-1}(2d)^{k}D(n-k-1)-\sum_{k=1}^{n-1}(2d)^{k}D(n-k-1)\Big)$\\[0.2cm]
& $=$ & $aD(n-1),$\\
\end{tabular}
$$
or equivalently
$$ D(n)=(2d+a)D(n-1).$$
But, this implies that $D(n)=a(2d+a)^{n-1}$. \end{proof}

\begin{corollary}\label{cor3}
Let $\alpha=(\alpha_i)_{i\geq 0}$ be a arithmetical sequence with
$\alpha_i=id$, and let $\beta=(\beta_{i})_{i\geq 0}$ be the square
sequence, i.e., $\beta_i=i^2$. If $D(n)=\det(P_{\alpha,
\beta}(n))$, then we have
\begin{equation}\label{new-6-6}
D(n)=-dD(n-2)+2d^2D(n-3).
\end{equation}
\end{corollary}
\begin{proof} By Corollary \ref{th2}, we get
$$D(n)=\sum_{k=0}^{n-1}(-d)^{k}\hat{\beta}_kD(n-k-1).$$
We claim that $\hat{\beta}_0=0$, $\hat{\beta}_1=1$,
$\hat{\beta}_2=2$ and $\hat{\beta}_k=0$ for $k\geq 3$. Now, it is
obvious that our claim implies the validity of Eq.
(\ref{new-6-6}).

Clearly $\hat{\beta}_0=0$, $\hat{\beta}_1=1$ and
$\hat{\beta}_2=2$. Now, we assume that $k\geq 3$. In this case we
have
$$\hat{\beta}_{k}=\sum_{j=0}^{k}(-1)^{k+j}{k\choose j}\beta_j=\sum_{j=0}^{k}(-1)^{k+j}{k\choose j}j^2.$$
We define the functions $f$ and $g$ as follows:
$$ f(x)=(1-x)^{k}=\sum_{j=0}^k{k \choose j}(-x)^{j} \ \ \ {\rm and} \ \ \ g(x)=-xf'(x).$$
Now, an easy calculation shows that
$$g'(x)=(1-x)^{k-2}\cdot(\ast)=\sum_{j=0}^{k}{k \choose j}j^2(-x)^{j-1},$$
and putting $x=1$, we get
$$\sum_{j=0}^{k}{k \choose j}j^2(-1)^{j-1}=0.$$
Now, by multiplying both sides by $(-1)^{k+1}$, we obtain
$$\sum_{j=0}^{k}{k\choose j}j^2(-1)^{k+j}=0,$$
as claimed. \end{proof}

\begin{corollary}\label{geometry} Let $\alpha=(\alpha_i)_{i\geq
0}$ and $\beta=(\beta_j)_{j\geq 0}$ be two geometric sequences
with $\alpha_i=\rho^{i}$ and $\beta_j=\sigma^j$. Then, we have
$$\det(T_{\alpha, \beta}(n))=(1-\rho\sigma)^{n-1}.$$
\end{corollary}
\begin{proof}
By Theorem \ref{th1}, we have $\det(T_{\alpha,
\beta}(n))=\det(P_{\check{\alpha}, \check{\beta}}(n))$. On the
other hand, straightforward computations show  that
$\check{\alpha}=(\check{\alpha}_i)_{i\geq 0}$ with
$\check{\alpha}_i=(1+\rho)^i$ and similarly
$\check{\beta}=(\check{\beta}_j)_{j\geq 0}$ with
$\check{\beta}_j=(1+\sigma)^j$. By applying Lemma 3, we conclude
the assertion.
\end{proof}

\subsection{Certain Generalized Pascal Triangles}

\begin{proposition}
Let $a, b, c\in \mathbb{C}$ and let $n$ be a positive integer. Let
$\alpha=(\alpha_i)_{i\geq 0}$ and $\beta=(\beta_j)_{j\geq 0}$ be
two sequences with $\alpha_i=(2^i-1)a+c$ and $\beta_j=(2^j-1)b+c$.
Then, we have
$$\det(P_{\alpha, \beta}(n))=\left\{\begin{array}{lll}
\lfloor \frac{1}{n}\rfloor c & \mbox{if} & a=b=c, \\[0.2cm]
\big[c+a(n-1)\big](c-a)^{n-1} & \mbox{if} & a=b\neq c, \\[0.2cm]
\frac{b}{b-a}(c-a)^{n}+\frac{a}{b-a}(c-b)^{n} & \mbox{if} & a\neq b. \\[0.2cm]
 \end{array} \right.$$
\end{proposition}
\begin{proof}
By Theorem \ref{th1}, we have $\det(P_{\alpha,
\beta}(n))=\det(T_{\hat{\alpha}, \hat{\beta}}(n))$. A
straightforward computation shows that
$$ \hat{\alpha}=(c, a, a, a, \ldots) \ \ \ \ {\rm and} \ \ \ \   \hat{\beta}=(c, b, b, b, \ldots).$$
Therefore, with notation in \cite{GST}, we have $T_{\hat{\alpha},
\hat{\beta}}(n)=M_n(b,a,c)$, and since $M_n(a,b,c)=M_n(b,a,c)^t$
we have
$$\det(T_{\hat{\alpha},
\hat{\beta}}(n))=\det(M_n(a,b,c)).$$ Now, the proof follows the
lines in the proof of Theorem 2 in \cite{GST}.
\end{proof}

\begin{proposition}
Let $a, b, c\in \mathbb{C}$ and let $n$ be a positive integer. Let
$\alpha=(\alpha_i)_{i\geq 0}$ and $\beta=(\beta_j)_{j\geq 0}$ be
two sequences with $\alpha_i=2^{i-1}(ia+2c)$ and
$\beta_j=2^{j-1}(jb+2c)$. Then, we have
$$\det(P_{\alpha, \beta}(n))=(-1)^{n+1}(a+b)^{n-2}\big[c(a+b)+(n-1)ab\big].$$
\end{proposition}
\begin{proof}
Again from Theorem \ref{th1}, we have $\det(P_{\alpha,
\beta}(n))=\det(T_{\hat{\alpha}, \hat{\beta}}(n))$, where
$\hat{\alpha}$ and  $\hat{\beta}$ are two arithmetical sequences
as
$$\hat{\alpha}=(\hat{\alpha}_i)_{i\geq 0}=(c, c+a, c+2a, \ldots, c+ia, \ldots) \ \ \ \ {\rm and} \ \ \ \
\hat{\beta}=(\hat{\beta}_j)_{j\geq 0}=(c, c+b, c+2b, \ldots,
c+jb, \ldots).  $$ Now we compute the determinant of
$T_{\hat{\alpha}, \hat{\beta}}(n)$. To do this, we apply the
following elementary column operations:
$${\rm C}_j\longrightarrow {\rm C}_j-{\rm C}_{j-1},
\ \ \ \ \ j=n-1, n-2,\ldots,2;$$ and we obtain the following
quasi-T$\rm \ddot{o}$eplitz matrix:
$$
\left(
\begin{array}{l|cccc}
c & b & b & \ldots & b \\
 \hline c+a &  &  & \\
c+2a &  &  &  & \\
\vdots &  &  &T_{\lambda, \mu}(n-1) &\\
c+(n-1)a &  &  & \\
\end{array} \right),$$
where $\lambda=(-a, -a, -a, \ldots)$ and $\mu=(-a, b, b,
\ldots)$.  Again, we subtract column $j$ from column $j+1$,
$j=n-2, n-3, \ldots, 2$. It is easy to see that, step by step,
the rows and columns are ``emptied" until finally the determinant
$$
\det(T_{\hat{\alpha}, \hat{\beta}}(n))=\det \left(
\begin{array}{lc|ccc}
c & b & 0 & \ldots & 0 \\
 \hline c+a & -a &  & \\
c+2a &  -a &  & (a+b)I_{(n-2)\times (n-2)} & \\
\vdots &  \vdots&  & &\\
c+(n-2)a & -a &  &  &  \\
\hline
c+(n-1)a & -a & 0 & \ldots & 0 \\
\end{array} \right),$$
is obtained. The proposition follows now immediately, by
expanding the determinant along the last row.
\end{proof}

\subsection{Fibonacci and Lucas Numbers as Principal Minors of a Quasi-Pascal Matrix
} There are several infinite matrices that the principal minors of
which form a Fibonacci or Lucas (sub)sequences. For instance, in
\cite{mpss}, we have presented a family of tridiagonal matrices
with the following form:
\begin{equation}\label{e0}
F_\lambda(\infty)=\left ( \begin{array}{cccccc} 1 & \lambda_0 & 0
& 0 & 0 &
\cdots\\
-\lambda_0^{-1} & 1 & \lambda_1 & 0 & 0 &
\cdots\\
0 & -\lambda_1^{-1} & 1 & \lambda_2 & 0 &
\cdots\\
0 & 0&-\lambda_2^{-1} &  1 & \ddots &
\cdots\\
\vdots & \vdots & \ddots & \ddots & \ddots &
\cdots\\
\end{array} \right )
\end{equation}
where $\lambda=(\lambda_i)_{i\geq 0}$ with $\lambda_i\in
\mathbb{C}^\ast=\mathbb{C}\backslash \{0\}$. Indeed, the principal
minors of these matrices for every $\lambda$ form the sequence
$\big(F(n+1)\big)_{n\geq 1}$ (Theorem 1 in \cite{mpss}). Also,
for the special cases $\lambda_0=\lambda_1=\ldots=\sqrt{\pm1}$
see \cite{Cahill1, Cahill2} and \cite{St1}. In (\cite{St2}, p.
555--557), Strang presents the infinite tridiagonal (T$\rm
\ddot{o}$eplitz) matrices:
\begin{equation}\label{1387-09-08}
P=T_{(3, t, 0, 0, \ldots), (3, t, 0, 0, \ldots)}(\infty),
\end{equation}
where $t=\pm 1$, and it is easy to show that the principal minors
of $T$ form the subsequence $\big(F(2n+2)\big)_{n\geq 1}$ from
Fibonacci sequence. As an another example, the principal minors of
T$\rm \ddot{o}$eplitz matrices:
\begin{equation}\label{1387-09-08}
Q=T_{(2, t, 1, 1, \ldots),(2, t, 0, 0, \ldots)}(\infty),
\end{equation}
where $t=\pm 1$, form the sequence $\big(F(n+2)\big)_{n\geq 1}$
for $t=1$ and the sequence $\big(F(2n+1)\big)_{n\geq 1}$ for
$t=-1$ (\cite{Cahill1}, Examples 1, 2).

We can summarize the above results in the following proposition.

\begin{proposition}\label{cahill-et-all}   (\cite{Cahill1, Cahill2, St1, St2})  Let $n$ be a natural number,
$\alpha=(\alpha_i)_{i\geq 0}$ and $\beta=(\beta_i)_{i\geq 0}$ be
two sequences, and let $d_n$ be the principal minor of
$T_{\alpha, \beta}(\infty)$. Then, the following hold.

$(1)$ If $\alpha=\beta=(1, \sqrt{-1}, 0, 0, \ldots )$, then
$d_n=F(n+1)$.

$(2)$ If $\alpha=\beta=(3, t, 0, 0, 0, \ldots)$ where $t=\pm 1$,
then $d_n=F(2n+2)$.

$(3)$ If $\alpha=(1, -1, 0, 0, \ldots)$ and $\beta=(1, 1, 0, 0,
\ldots)$, then $d_n=F(n+1)$.

$(4)$ If $\alpha=(2, 1, 1, 1, \ldots)$ and $\beta=(2, -1, 0, 0,
\ldots)$, then $d_n=F(2n+1)$.

$(5)$ If $\alpha=(2, 1, 1, 1, \ldots)$ and $\beta=(2, 1, 0, 0,
\ldots)$, then $d_n=F(n+2)$.
\end{proposition}

Let $\phi=\frac{1+\sqrt{5}}{2}$, the golden ratio, and
$\Phi=\frac{1-\sqrt{5}}{2}$, the golden ratio conjugate. The
recent paper of Griffin, Stuart and Tsatsomeros \cite{GST} gives
the following result:
\begin{proposition}\label{GriffinST}   (\cite{GST}, Lemma 7)
For each positive integer $n$, let
$$P_n=T_{(1, \Phi, \Phi, \ldots), (1, \phi, \phi, \ldots)}(n),
\ \ {\rm and} \ \ Q_n=T_{(0, -\Phi, -\Phi, \ldots), (0, -\phi,
-\phi, \ldots)}(n).$$ Then we have
$$ \det(P_n)=F(n+1), \ \  {\rm and} \ \ \det(Q_n)=F(n-1).$$
\end{proposition}

Using Propositions \ref{cahill-et-all}, \ref{GriffinST} and
Theorem \ref{th1}, we immediately deduce the following corollary.

\begin{corollary}\label{generalized-Fibonacci} Let $n$ be a natural number,
$\alpha=(\alpha_i)_{i\geq 0}$ and $\beta=(\beta_i)_{i\geq 0}$ be
two sequences, and let $d_n$ be the principal minor of $P_{\alpha,
\beta}(\infty)$. Then, the following hold.

$(1)$ If $\alpha_i=\beta_i=1+i\sqrt{-1}$, then $d_n=F(n+1)$.

$(2)$ If $\alpha_i=\beta_i=3-i$, then $d_n=F(2n+2)$.

$(3)$ If $\alpha_i=\beta_i=3+i$, then $d_n=F(2n+2)$.

$(4)$ If $\alpha_i=1-i$ and $\beta_i=1+i$, then $d_n=F(n+1)$.

$(5)$ If $\alpha_i=2^i+1$ and $\beta_i=2-i$, then $d_n=F(2n+1)$.

$(6)$ If $\alpha_i=2^i+1$ and $\beta_i=2+i$, then $d_n=F(n+2)$.

$(7)$ If $\alpha_i=(2^i-1)\Phi+1$ and $\beta_i=(2^i-1)\phi+1$,
then $d_n=F(n+1)$.

$(8)$ If $\alpha_i=(1-2^i)\Phi$ and $\beta_i=(1-2^i)\phi$, then
$d_n=F(n-1)$.
\end{corollary}

In the sequel, we study together the sequences $\mathcal{F}$ and
$\mathcal{L}$, and, in order to unify our treatment, we introduce
the following useful notations. For $\varepsilon \in \{+, -\}$ we
let $F^{\varepsilon}(n)=F(n)$ if $\epsilon =+$; and
$F^{\varepsilon}(n)=L(n)$ if $\epsilon =-$.

\begin{theorem}\label{th1387-9-21}
Let $r$ be a non-negative integer and $s$ be a positive integer.
Suppose that
$$\phi_{r,s}=\Big\lceil\frac{F^{\varepsilon}(2r+s)}{F^{\varepsilon}(r+s)}\Big\rceil \ \ \ \mbox{and} \ \ \  \psi_{r,s}=
\sqrt{\phi_{r,s}F^{\varepsilon}(r+s)-F^{\varepsilon}(2r+s)}.$$
Then, the principal minors of the following infinite quasi-Pascal
matrix:
$$P^{[r,s]}(\infty)=\left [\begin{array}{cc|ccc}
F^{\varepsilon}(r+s) & \psi_{r,s} & 0 & 0 & \ldots  \\
\psi_{r,s} & \phi_{r,s} & \sqrt{(-1)^r} & \sqrt{(-1)^r} & \ldots
\\ \hline &&&& \\[-0.4cm]
0 & \sqrt{(-1)^r} &  & &  \\
0 & \sqrt{(-1)^r} & & P_{\alpha, \alpha}(\infty) & \\
\vdots  & \vdots & &  &  \\
\end{array} \right]$$
where $\alpha=(\alpha_i)_{i\geq 0}$ is an arithmetical sequence
with $\alpha_i=F^{-}(r)+i\sqrt{(-1)^r}$, form the subsequence
$\{F^{\varepsilon}(nr+s)\}_{n=1}^{\infty}$ from Fibonacci or
Lucas sequences.
\end{theorem}
\begin{proof} Cahill and Narayan in \cite{Cahill3} introduce the
following quasi-T$\rm \ddot{o}$eplitz matrices:

$$T^{[r,s]}(\infty)=\left [\begin{array}{cc|ccc}
F^{\varepsilon}(r+s) & \psi_{r,s} & 0 & 0 & \ldots \\  \psi_{r,s}
& \phi_{r,s} & \sqrt{(-1)^r} & 0 & \ldots
\\ \hline &&&& \\[-0.35cm]
0 & \sqrt{(-1)^r} &  & &  \\
0 & 0 & & T_{\beta, \beta}(\infty) & \\
\vdots  & \vdots & &  &  \\
\end{array} \right]$$
where $\beta=(F^{-}(r), \sqrt{(-1)^r}, 0, 0, \ldots )$. Moreover,
they show that
\begin{equation}\label{cahill332}\det(T^{[r,s]}(n))=F^{\varepsilon}(nr+s).\end{equation}
Now, we decompose the matrix $T^{[r,s]}(n)$ as follows:
\begin{equation}\label{cahill333}
T^{[r,s]}(n)=\tilde{L}(n)P^{[r,s]}(n)\tilde{L}(n)^t,
\end{equation} where
$$\tilde{L}(n)=I_{2\times 2}\oplus L^{-1}(n-2).$$
The proof of Eq. (\ref{cahill333}) is similar to the proof of
Theorem \ref{th1} and we omit it here. Now, using Eqs.
(\ref{cahill332}) and (\ref{cahill333}), we easily see that
$$\det(P^{[r,s]}(n))=F^{\varepsilon}(nr+s),$$
and the proof of theorem is complete.
\end{proof}

\subsection{Generalized Pascal Triangle Associated to a Constant Sequence}
Another consequence of Theorem \ref{th1}, is the following.
\begin{corollary}\label{cor1}
Let $\alpha=(\alpha_i)_{i\geq 0}$ and $\beta=(\beta_i)_{i\geq 0}$
be two sequences with $\alpha_0=\beta_0=\gamma$. If $\alpha$ or
$\beta$ is a constant sequence, then we have
$\det(P_{\alpha,\beta}(n))=\gamma^n$.
\end{corollary}
\begin{proof} By Theorem \ref{th1}, we have $\det(P_{\alpha,\beta}(n))=\det(T_{\hat{\alpha},\hat{\beta}}(n))$.
But in both cases, the T$\rm \ddot{o}$eplitz matrix
$T_{\hat{\alpha},\hat{\beta}}(n)$ is a lower triangular matrix or
an upper triangular one with $\gamma$ on its diagonal. This
implies the corollary.
\end{proof}

The generalized Pascal triangle $P_{\alpha,\alpha}(\infty)$
associated to the pair of identical sequences $\alpha$ and
$\alpha$, is called the {\em generalized symmetric Pascal
triangle} associated to $\alpha$ and yields symmetric matrices
$P_{\alpha,\alpha}(n)$ by considering principal submatrices
consisting of the first $n$ rows and columns of
$P_{\alpha,\alpha}(\infty)$.  For an arbitrary sequence
$\alpha=(\alpha_i)_{i\geq 0}$ with $\alpha_0=0$, we define
$\tilde{\alpha}=(\tilde{\alpha}_i)_{i\geq 0}$ where
$\tilde{\alpha}_i=(-1)^{i}\alpha_i$ for all $i$. Then, the
generalized Pascal triangle $P_{\alpha, \tilde{\alpha}}(\infty)$
associated to the sequences $\alpha$ and $\tilde{\alpha}$, is
called the {\em generalized skymmetric Pascal triangle} associated
to $\alpha$ and $\tilde{\alpha}$, and yields skymmetric matrices
$P_{\alpha, \tilde{\alpha}}(n)$ by considering principal
submatrices consisting of the first $n$ rows and columns of $P_{
\alpha, \tilde{\alpha}}(\infty)$.

\begin{ex}\label{ex1}
Let $n\geq 2$ be a natural number. Then

$(i)$ The generalized symmetric Pascal triangle $P_{\mathcal{F},
\mathcal{F}}(n)$ has determinant $-2^{n-2}$.

$(ii)$ The generalized skymmetric Pascal triangle
$P_{\mathcal{F},\tilde{\mathcal{F}}}(n)$ has determinant
$2^{n-2}$.
\end{ex}

All assertions in this Example follow of course from Theorem 3.1
in \cite{bacher}. However, we will reprove them independently.

\begin{proof} $(i)$ Consider the generalized symmetric Pascal triangle
$$P_{\mathcal{F},
\mathcal{F}}(n)=\left( \begin{array}{ccccccc}
0 & 1 & 1 & 2& 3&   \ldots & \mathcal{F}_{n-1} \\
1 & 2 & 3 & 5 & 8  &   \ldots & \mathcal{F}_{n+1}  \\
1 & 3 & 6 & 11& 19 &   \ldots & . \\
2 & 5 & 11 & 22& 41 &   \ldots & . \\
3 & 8 & 19 & 41& 82 &   \ldots & . \\
\vdots & \vdots & \vdots & \vdots & \vdots &   \ddots & . \\
\mathcal{F}_{n-1} & \mathcal{F}_{n-2} & . &  . & . &   \ldots & . \\
\end{array} \right)$$
Now, we apply the following elementary row operations: $${\rm
R}_i\longrightarrow {\rm R}_i-{\rm R}_{i-1}-{\rm R}_{i-2}, \ \ \
\ \ i=n-1, n-2, \ldots, 2.$$ It is easy to see that
$$
\begin{tabular}{lll}
$\det\big(P_{\mathcal{F}, \mathcal{F}}(n)\big)$& $=$ & $\det\left(
\begin{array}{cc|ccccc}
0 & 1 & 1 & 2& 3&   \ldots & \mathcal{F}_{n-1} \\
1 & 2 & 3 & 5 & 8  &   \ldots & \mathcal{F}_{n+1}\\ \hline
0 & 0 & 2 & 4& 8 &   \ldots & 2(\mathcal{F}_{n}-1)\\
0 & 0 & 2 & 6& 14 &   \ldots & \ast\\
0 & 0 & 2 & 8& 22 &   \ldots & \ast\\
\vdots & \vdots & \vdots & \vdots & \vdots &  \ddots & \ast \\
0 & 0 & 2& \ast&\ast &\ldots&\ast \\
\end{array} \right)$ \\[2cm]
&$ =$ & $\det\left( \begin{array}{cc|ccccc}
0 & 1 & 1 & 2& 3&   \ldots & \mathcal{F}_{n-1}\\
1 & 2 & 3 & 5 & 8  &   \ldots & \mathcal{F}_{n+1}\\ \hline
0 & 0 &  & &   &   \\
0 & 0 &  &   & &   \\
\vdots & \vdots   & &  &  & P_{\lambda, \mu}(n-2)& \\
0 & 0 & & & & & \\
\end{array} \right)$\\
\end{tabular}$$ where $\lambda=(2, 2, 2, \ldots)$ and $\mu=(2(\mathcal{F}_3-1), 2(\mathcal{F}_4-1), 2(\mathcal{F}_5-1),
2(\mathcal{F}_6-1), \ldots)$. Now, by Corollary \ref{cor1}, we get
$$\det(P_{\mathcal{F}, \mathcal{F}}(n))=\det\left(\begin{array}{cc}0&1\\1&2\end{array}\right)
\cdot \det\big(P_{\lambda, \mu}(n-2)\big)=-2^{n-2},$$ as desired.

$(ii)$ Here, we consider the generalized skymmetric Pascal
triangle
$$P_{\mathcal{F},
\tilde{\mathcal{F}}}(n)=\left( \begin{array}{ccccccc}
0 & -1 & 1 & -2& 3&   \ldots & (-1)^{n-1}\mathcal{F}_{n-1} \\
1 & 0 & 1 & -1 & 2  &   \ldots &  . \\
1 & 1 & 2 & 1& 3 &   \ldots &  . \\
2 & 3 & 5 & 6& 9 &   \ldots &  . \\
3 & 6 & 11 & 17& 26 &   \ldots&  . \\
\vdots & \vdots & \vdots & \vdots &  \vdots &  \ddots &  . \\
\mathcal{F}_{n-1}  & . & . &  . & .  &  \ldots &  . \\
\end{array} \right)$$
Similarly, we apply the following elementary column operations:
$${\rm C}_j\longrightarrow {\rm C}_j+{\rm C}_{j-1}-{\rm C}_{j-2},
\ \ \ \ \ j=n-1,n-2,\ldots,2;$$ and we obtain
$$
\begin{tabular}{lll}
$\det(P_{\mathcal{F}, \tilde{\mathcal{F}}}(n))$& $=$ & $\det\left(
\begin{array}{cc|ccccc}
0 & -1 & 0 & 0& 0&   \ldots & 0 \\
1 & 0 & 0 & 0 & 0  &  \ldots & 0 \\
\hline
1 & 1 & 2 & 2& 2 &   \ldots & 2\\
2 & 3 & 6 & 8& 10 &   \ldots & \ast\\
3 & 6 & 14 & 22& 32 &   \ldots & \ast\\
\vdots & \vdots & \vdots &  \vdots & \vdots &  \ddots & \ast \\
\mathcal{F}_{n-1} & \mathcal{F}_{n+1}-2 & \ast& \ast &\ast&\ldots&\ast \\
\end{array} \right)$ \\[2cm]
&$ =$ & $\det\left( \begin{array}{cc|ccccc}
0 & -1 & 0 & 0& 0&   \ldots & 0\\
1 & 0 & 0 & 0 & 0  &   \ldots & 0\\
\hline
1 & 1 &  & &   &   \\
2 & 3 &  &   & &   \\
\vdots & \vdots   & &  &  & P_{\nu, \lambda}(n-2)& \\
\mathcal{F}_{n-1} & \mathcal{F}_{n+1}-2 & & & & & \\
\end{array} \right)$\\
\end{tabular}$$ where $\lambda=(2, 2, 2, \ldots)$. Again, by Corollary \ref{cor1}, we get
$$\det(P_{\mathcal{F}, \tilde{\mathcal{F}}}(n))=\det\left(\begin{array}{cc}0&1\\-1&0\end{array}\right)
\cdot \det\big(P_{\nu,\lambda}(n-2)\big)=2^{n-2},$$ as desired.
\end{proof}

\begin{ex}\label{ex2}
Let $n\geq 2$ be a natural number. Then the generalized Pascal
triangle $P_{\mathcal{F}^\ast, \mathcal{I}^\ast}(n)$ has
determinant $(-1)^n$.
\end{ex}
\begin{proof}
Consider the following generalized Pascal triangle:
$$P_{\mathcal{F}^\ast, \mathcal{I}^\ast}(n)=\left( \begin{array}{ccccccc}
1 & 2 & 6 & 24& 120&   \ldots & n! \\
1 & 3 & 9 & 33 & 153  &   \ldots & . \\
2 & 5 & 14 & 47& 200 &   \ldots & . \\
3 & 8 & 22 & 69 & 269 &   \ldots & . \\
5 & 13 & 35 & 104 & 373 &   \ldots & . \\
\vdots & \vdots & \vdots & \vdots & \vdots &   \ddots & . \\
\mathcal{F}_{n-1} & \mathcal{F}_{n+1} & . & . & .&  \ldots & . \\
\end{array} \right)$$
Again, we use the similar elementary row operations as Example
$1(i)$:
$${\rm R}_i\longrightarrow {\rm R}_i-{\rm R}_{i-1}-{\rm R}_{i-2},
\ \ \ \ \ i=n-1, n-2, \ldots, 2.$$ Therefore, we deduce that
$$
\begin{tabular}{lll}
$\det\big(P_{\mathcal{F}^\ast, \mathcal{I}^\ast}(n)\big)$& $=$ &
$\det\left(
\begin{array}{cc|ccccc}
1 & 2 & 6 & 24& 120 & \ldots & n! \\
1 & 3 & 9 & 33 & 153  &  \ldots & \ast\\ \hline
0 & 0 & -1 & -10 & -73 &  \ldots & \ast\\
0 & 0 & -1 & -11 & -84 &  \ldots & \ast\\
0 & 0 & -1 & -12 & -96 & \ldots & \ast\\
0 & 0 & -1 & -13 & -109&  \ldots & \ast\\
\vdots & \vdots & \vdots & \vdots & \vdots &  \ddots & \ast \\
0 & 0 & -1& \ast&\ast &\ldots&\ast \\
\end{array} \right)$ \\[2cm]
&$ =$ & $\det\left( \begin{array}{cc|ccccc}
1 & 2 & 6 & 24& 120 &   \ldots & n!\\
1 & 3 & 9 & 33 & 153  &   \ldots & \ast\\ \hline
0 & 0 &  & &   &   \\
0 & 0 &  &   & &   \\
\vdots & \vdots   & &  &  & P_{\lambda, \mu}(n-2)& \\
0 & 0 & & & & & \\
\end{array} \right)$\\
\end{tabular}$$ where $\lambda=(-1, -1, -1, \ldots)$ and $\mu=(-1, -10, -73,
\ldots)$. Now, by Corollary \ref{cor1}, we get
$$\det\big(P_{\alpha, \beta}(n)\big)=\det\left(\begin{array}{cc}1&2\\1&3\end{array}\right)
\cdot \det\big(P_{\lambda, \mu}(n-2)\big)=(-1)^{n-2}=(-1)^n,$$ as
desired.
\end{proof}

{\bf Acknowledgment.} The first author would like to thank IPM
for the financial support.


\begin{thebibliography}{99}
\bibitem{bacher} R. Bacher. Determinants of matrices related to
the Pascal triangle. {\em J. Theorie Nombres Bordeaux,} 14:19-41,
2002.

\bibitem{BC} R. Bacher and R. Chapman. Symmetric Pascal
matrices modulo $p$. {\em European J. Combinatorics,} 25:459-473,
2004.

\bibitem{Cahill1} N. D. Cahill, J. R. D'Errico, D. A. Narayan, and J. Y.
Narayan. Fibonacci determinants, {\em The College Math. J.},
33(3)(2002), 221-225.

\bibitem{Cahill2} N. D. Cahill, J. R. D'Errico and J. P. Spence.
Complex factorizations of the Fibonacci and Lucas numbers. {\em
Fibonacci Quart.}, 41(1)(2003), 13-19.

\bibitem{Cahill3} N. D. Cahill and D. A. Narayan.
Fibonacci and Lucas numbers as tridiagonal matrix determinants.
{\em Fibonacci Quart.}, 42(3)(2004), 216-221.


\bibitem{GST} K. Griffin, J. L. Stuart and M. J. Tsatsomeros.
Noncirculant T$\rm \ddot{\bf o}$eplitz matrices all of whose
powers are T$\rm \ddot{\bf o}$eplitz, {\em Czechoslovak
Mathematical Journal}, 58(4) 1185-1193, 2008.

\bibitem{IS} Integer-sequences,
http://www.research.att.com/~njas/sequences/index.html

\bibitem{kratt} C. Krattenthaler. Advanced determinant
calculus. {\em  S$\acute{e}$minaire Lotharingien Combin.,}
Article B42q, 67 pp., (1999).

\bibitem{kkk} C. Krattenthaler. Evaluations of some
 determinants of matrices related to the
Pascal triangle, {\em Semin. Lothar. Comb.,} Article B47g, 19 pp,
2002.

\bibitem{edelman} A. Edelman and G. Strong. Pascal matrices.
{\em Amer. Math. Monthly,} 111:189-197, 2004.

\bibitem{mpss} A. R. Moghaddamfar, S. M. H.
Pooya, S. Navid Salehy and S. Nima Salehy. Fibonacci and Lucas
sequences as the principal minors of some infinite matrices.
Submitted for publication.


\bibitem{St1} G. Strang. Introduction to Linear Algebra, Third Edition.
Wellesley-Cambridge Press, 2003.

\bibitem{St2} G. Strang and K. Borre. Linear Algebra, Geodesy, and
GPS. Wellesley-Cambridge Press, 1997.
\end{thebibliography}
\end{document}